# Basic Dynamical Systems Control of Aeroengine Flow

Björn Birnir[1] and Höskuldur A. Hauksson

Department of Mathematics, University of California Santa Barbara, Santa Barbara, CA 93106
`birnir@math.ucsb.edu, http://math.ucsb.edu/~birnir`

**Abstract**

The notions of basic controllability and basic control are defined using dynamical systems theory of partial differential equations. A quadratic optimal control of the linearized viscous Moore-Greitzer equation is presented and it is confirmed that stall is uncontrollable in this model. A basic control is constructed for the nonlinear viscous Moore-Greitzer equation which can control both surge and stall. Numerical simulations of the basic control are presented.

## 1 Introduction

In recent years a lot of attention has been devoted to the study of air flow through turbomachines. The main reason for this interest is that when a turbomachine, such as a jet engine, operates close to its optimal operating parameter values, the flow can become unstable. These instabilities put a large stress on the engine and in some cases the engine needs to be turned off in order to recover original operation. For this reason jet engines are currently operated away from their optimal operating parameter values.

Moore and Greitzer published in 1986 a PDE model for turbomachines which has been very successful [17]. A substantial amount of work has been done on finite Galerkin approximations of that model since, see e.g. [15] [13] and references therein. Banaszuk *et al.* considered the full PDE model of Moore and Greitzer. The model is known as the viscous Moore-Greitzer equation (vMG equation), when a viscous term is added.

Birnir and Hauksson [5] proved that the vMG equation is well posed in the Hilbert space $X = \bar{H}^1 \times \mathbf{R}^2$ where $\bar{H}^1$ denotes the Sobolev space with index one of functions on the unit circle with square integrable first derivative and zero mean. This solution is smooth in space and time variables and this dynamical system has a global attractor with finite Hausdorff and fractal dimensions. In [4] the authors analyzed the basic attractor and found explicit solutions for stall for certain parameter values and showed that they are stable and persist under small perturbations of the parameters. Stall is a solitary wave that rotates around the annulus at half the rotor speed of the engine. They conclude that the basic attractor consists of design flow, surge and one or more stall solutions. The analysis of the basic attractor (see Subsection 1.1) was extended for all parameter values in [6] and there they derived a reduced order model that captures the dynamics of the vMG equation quantitatively as well as qualitatively. These results are in good agreement with experimental and numerical results [10]. In [7] Birnir and Hauksson addressed the controllability of design flow, stall and surge.

The backstepping control given by Banaszuk *et al.*[1] shows that one can eliminate stall and surge by using throttle control. The question we want to answer in this paper is how simple can we make the control design and how efficient can the control be? The control philosophy we want to adopt, is to construct a control strategy that can recover design flow operation after large disturbances, but is not necessarily good for regulating the design flow. For that a different strategy would be used.

In [7] we prove that the vMG equation with throttle control is not basically controllable. Hoever, if one in addition to throttle control has air injection or bleeding at ones disposal, then the vMG equation is basically controllable. More details can be found in [7].

### 1.1 The Equation of Motion and Assumptions

Moore and Greitzer [17] derived a model of the three dimensional flow through the compression system in jet engines. When one assumes that the flow does not depend on the radial direction, the equations reduce to the following.

$$\frac{\partial}{\partial t}\varphi = \nu\frac{\partial^2\varphi}{\partial\theta^2} - \frac{1}{2}\frac{\partial\varphi}{\partial\theta} + \psi_c(\Phi+\varphi) - \overline{\psi_c}, \ \theta \in [0, 2\pi) \quad (1)$$

$$\dot{\Phi} = \frac{1}{l_c}(\overline{\psi_c} - \Psi) \quad (2)$$

$$\dot{\Psi} = \frac{1}{4l_c B^2}(\Phi - \gamma F_T^{-1}(\Psi)) \quad (3)$$

$$\quad (4)$$

where

$$\overline{\psi_c} := \frac{1}{2\pi}\int_0^{2\pi} \psi_c(\Phi(t) + \varphi(t,\theta))d\theta.$$

This equation is known as the viscous Moore-Greitzer equation (vMG equation). Here the dot represents the total derivative with respect to time.

The characteristic $\psi_c$ is a cubic polynomial with a negative leading coefficient and $F_T^{-1}$ is a smooth function which is equal to $F_T^{-1}(\Psi) = \Psi|\Psi|^{-1/2}$ outside a small neighbourhood of the origin.

---

[1] and Science Institute, University of Iceland, 3 Dunhaga, IS-107 Reykjavik, Iceland

In the sequel we will allow γ to depend on the state but we will assume that it does so in a smooth way and that there exists a constant $\tilde{\gamma}$ such that

$$\gamma(\Phi, \Psi, \varphi) \geq \tilde{\gamma} > 0.$$

With these restrictions on γ the results on existence of unique solutions and their regularity [5] still hold. In addition, the system will again have a global attractor whose fractal and Hausdorff dimensions can be bounded by the same bounds as in [5] with γ replaced by $\tilde{\gamma}$. The global attractor consists of a finite-dimensional set of solutions that all other solutions tend towards as time progresses.

## 2 The Basic Attractor

Once the existence of a global attractor has been established, the natural question arises: How can one construct the global attractor and can one obtain a system of ODEs that describe the evolution on the attractor? There are, for the most part, two main approaches that researchers have taken here.

The first one, and the more popular one, is to think of the attractor as a set embedded in a larger manifold, often called intertial manifold, see e.g. [9]. The problem of finding ODEs that describe the flow on this manifold or an approximate manifold, is then solved by using a Galerkin projection onto a basis. The number of basis vectors needed is often quite large. This can be due to either that the bounds on the dimension of the attractors found by current methods tend to be rather conservative, or that the asymptotic dynamics of the system in question are in fact high dimensional. Hence the system of ODEs is not tractable for analytical analysis, but may lend itself better to numerical work.

The second approach is to consider only the core of the attractor called the basic attractor (see below). Here one constructs the particular solutions in the attractor which attract "almost all" of the phase space. For some systems, the asymptotic dynamics in this "almost every" sense are low dimensional and one can completely determine the flow on the basic attractor *analytically* or by using the qualitative analysis of ODE dynamical systems.

Here we adopt the second approach, but before we go further, let us clarify what we mean by the basic attractor and by "almost every".

### 2.1 Prevalence and Basic Attractors

We need to extend the measure theoretic terms measure zero and almost every to infinite dimensional Banach spaces. Furthermore, we want to do it in such a way that these definitions behave well under the operations of the vector space. It turns out that it suffices that they behave well under translations of the set. The problem here lies in that there do not exist any nontrivial translation invariant measures in infinite dimensional spaces. If a subset $U \subset X$ in an infinite dimensional Banach space is nonempty and $\mu$ is a translation invariant measure on $X$, then either $\mu(U) = 0$ or $\mu(U) = \infty$. Following Hunt *et al.* [12] the ideas of measure zero and almost every can be replaced by shy and prevalent.

**Definition 2.1** Let $X$ denote a separable Banach space. We denote by $S + v$ the translate of the set $S \subset X$ by a vector $v$. A measure $\mu$ is said to be *transverse* to a Borel set $S \subset X$ if the following two conditions hold:

- There exists a compact set $U \subset X$ for which $0 < \mu(U) < \infty$.
- $\mu(S + v) = 0$ for every $v \in X$.

A Borel set $S \subset X$ is called *shy* if there exists a compactly supported measure transverse to $S$. More generally, a subset of $X$ is called shy if it is contained in a shy Borel set. The complement of a shy set is said to be a *prevalent* set.

The basic attractor should be the smallest part of the global attractor $\mathcal{A}$ which attracts a prevalent set. Let us make this more precise.

**Definition 2.2** An attractor $\mathcal{B}$ is a *basic attractor* if it satisfies the two conditions

1. The basin of attraction of $\mathcal{B}$ is prevalent.
2. $\mathcal{B}$ is minimal with respect to property (1), i.e. there exists no strictly smaller $\mathcal{B}' \subset \mathcal{B}$ with basin($\mathcal{B}$) $\subset$ basin($\mathcal{B}'$), up to shy sets.

This means that every point of $\mathcal{B}$ is essential, no point can be removed without removing a portion of the basin that is not shy. In numerical simulations or in physical experiments one would therefore only expect to observe the basic attractor after a long enough settling period.

In general, the basic attractor will be disconnected although the global attractor is connected. We can therefore speak of components of the basic attractor.

The following theorem, which is an extension of a finite dimensional version by Milnor [16], was proven in Birnir [2].

**Theorem 2.1** *Let $\mathcal{A}$ be the compact attractor of a continuous map $T(t)$ on a separable Banach space $X$. Then $\mathcal{A}$ can be decomposed into a basic attractor $\mathcal{B}$ and a remainder $\mathcal{C}$,*

$$\mathcal{A} = \mathcal{B} \cup \mathcal{C}$$

*such that basin($\mathcal{B}$) is prevalent and basin($\mathcal{C}$)\basin($\mathcal{B}$) is shy.*

It turns out that in the cases where an explicit description of $\mathcal{B}$ has been given, see [3] and [8], that the dimension of $\mathcal{B}$ is small, whereas the dimension of $\mathcal{C}$ can be quite large.

## 2.2 The Geometry of the Basic Attractor

Experimental and numerical evidence indicate that the basic attractor in axial compression systems is low dimensional [11] [18]. It can consist of a combination of axisymmetric design flow, surge and stall. The design flow is a stationary solution and surge is a periodic cycle which has been well studied (see [17] and [15]). It only involves the two ODEs in the system (1-3). Stall has been studied in a low order Galerkin truncations of the Moore-Greitzer equations [15] [13].

Under normal conditions the engine operates in design flow. There the flow through the compressor is uniform in space and time and the pressure rise is relatively high. In particular, $\varphi = 0$ and $\Phi$ and $\Psi$ are constant.

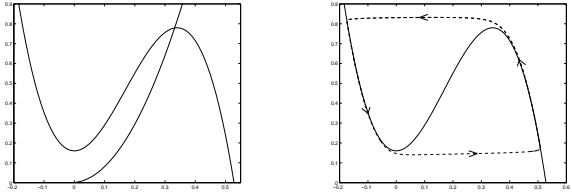

**Figure 1:** Left: Two characteristics: the cubic compressor characteristic and the parabolic throttle characteristic. Right: The surge limit cycle in the $(\Phi, \Psi)$ plane where $\varphi = 0$.

Figure 1 shows the $(\Phi, \Psi)$ plane. The parabola starting at the origin represents all stationary solutions for equation (3) and is called the throttle characteristic. The cubic curve represents all stationary solutions for equation (2), given that $\varphi = 0$, and is called the compressor characteristic. Since $\varphi = 0$ is a stationary solution for equation (1) we can conclude that the intersection of the two curves in Figure 3.1 is a stationary solution for the full system (1-3). This stationary solution is called design flow.

Design flow is stable to the right of the peak of the compressor characteristic. It is desirable to operate the engine on the right side of the peak with as high a pressure rise as possible without risking the system being thrown over to the unstable side by disturbances.

Surge is a limit cycle in the two ODEs (2-3) where the non-axisymmetric disturbance is zero, $\varphi = 0$. It has been studied by many authors, among them Greitzer [11] and Mc-Caughan [14] [15]. It arises as a subcritical Hopf bifurcation in the system

$$\dot{\Phi} = \frac{1}{l_c}(\psi_c(\Phi) - \Psi)$$
$$\dot{\Psi} = \frac{1}{4B^2 l_c}(\Phi - \gamma F_T^{-1}(\Psi))$$

which occurs for a large enough $B$ when the throttle parameter $\gamma$ is decreased. Since the bifurcation is subcritical, we have a one parameter family of unstable surge cycles that originate from the bifurcation point. This branch bends on itself and the cycles become stable [15]. These stable cycles are fairly large and a simulation of one is shown in Figure 1.

The solution spends most of its time on the two vertical sides of the cycle. There the slope of the compressor characteristic is negative so all nonaxisymmetric disturbances are damped.

Stall is a solitary wave solution. The wave rotates around the unit circle and the average flow $\Phi$ and pressure rise $\Psi$ are constant. When one looks for traveling wave solutions of the vMG equation, the problem can be reduced to finding periodic solutions of the Duffing's equation with the correct periods [4] [6]. These periodic solutions lie inside a homoclinic (or heteroclinic) orbit, and since the compressor characteristic is a cubic polynomial, these solutions can be found explicitly with quadratures. They can be expressed as rational functions of the Jacobi elliptic function $ns$ and varying the parameters in the vMG equation we can in fact construct a one parameter family of stall solutions, see [6].

## 3 Basic Controllability

Let us consider now the issue of controllability. In finite dimensional control theory, a system is said to be controllable if for every two points $x_0, x_1 \in X$ and every two real numbers $t_0 < t_1$, there exists a control function $u$ such that the unique solution of the equation

$$\dot{x} = F(x, u), \quad x(t_0) = x_0 \qquad (5)$$

satisfies $x(t_1) = x_1$.

In infinite dimensional spaces this notion of controllability is too restrictive. For practical control applications one can never have more than finitely many control parameters, if for no other reason, the fundamentals of computing require computer outputs to be finite. There is therefore no hope that nonlinear evolution equations in infinite dimensional spaces will be controllable in this strict sense in practical applications.

If an evolution equation has an attractor and a basic attractor its solutions will converge asymptotically to the attractor for all initial conditions and to the basic attractor for almost all initial conditions. The simplest thing one could ask of the control is that it make all or almost all initial conditions give rise to solutions that converge to a given component in the basic attractor. A more stringent requirement on the control would be that it make all or a prevalent set of (almost all) initial conditions give rise to solutions that converge to a given component in the global attractor. This requires one to have enough control authority over the local unstable manifolds of the hyperbolic trajectories in the attractor to make them attractive. Consider the following definitions

**Definition 3.1** The equation (5) is *basically controllable* if for every bounded set $M$ and every $\varepsilon > 0$, there exists a finite time $T(M)$ and a control function $u(t)$, such that for every solution $x(t)$ with initial data $x_o \in M$ and any minimal component of the basic attractor $\mathcal{B}_j$

$$\|x(t) - \mathcal{B}_j\| < \varepsilon,$$

for $t > T(M)$.

This definition say that given an initial point one can stear to any component of the basic attractor in finite time. It is hopeless to get a finite $T$ for $x_o$ lying in a prevalent (full measure) set in the infinite-dimensional space, for the reason discussed above. It is not wise to attempt to control every solution in the $\mathcal{A}$-attractor, because in general it ($\mathcal{C}$) contains many hyperbolic solutions and their heteroclinic connections.

**Definition 3.2** The equation (5) is *attractively controllable* if for every bounded set $M$ and every $\varepsilon > 0$, there exists a finite time $T(M)$ and a control function $u(t)$, such that for every solution $x(t)$ with initial data $x_o \in M$ and any trajectory $z$ in the attractor $\mathcal{A}$

$$\|x(t) - \omega(z)\| < \varepsilon,$$

for $t > T(M)$.

This definition says that given an initial point one can stear to the $\omega$ limit set of any trajectory in the $\mathcal{A}$-attractor in finite time.

One can also speak of basic controllability as b-controllability and attractive controllability as a-controllability. Clearly a-controllability implies b-controllability. Not surprisingly the control construction relies heavily on the geometry of the basic attractor. Consequently it is referred to as *basic control*. The remainder $\mathcal{C} = \mathcal{A}\setminus\mathcal{B}$ from Section 2.1 plays a large role in basic control. In general one would like to use its heteroclinic connections to move efficiently from one minimal basic attractor to another.

## 4 Basic Control for Design Flow of the Linearized Equation

The most important component of the basic attractor of the vMG equation is the design flow component. The goal is to construct a basic control that makes all solutions converge to the design flow. The simplest approach one could take would be to linearize the system about design flow, $(\Phi_0, \Psi_0, 0)$, corresponding to a throttle parameter $\gamma_0$, and apply the classical optimal control theory. We define the control parameter $u = \gamma - \gamma_0$ and we make a change of coordinates $(t, \eta) = (t, \theta - \frac{1}{2}t)$ to simplify the equations. Furthermore, we define the variable $y = (y_1, y_2, y_3) = (\varphi, \Phi - \Phi_0, \Psi - \Psi_0)$. The linearized equations can now be written as

$$\dot{y} = Ay + Bu \tag{6}$$

where

$$A = \begin{bmatrix} \nu \partial_\theta^2 + \psi'_c(\Phi_0) & 0 & 0 \\ 0 & \frac{1}{l_c}\psi'(\Phi_0) & -\frac{1}{l_c} \\ 0 & \frac{1}{4l_c B^2} & \frac{-\gamma_0}{4l_c B^2}(F_T^{-1})'(\Psi_0) \end{bmatrix} \tag{7}$$

and

$$B = \begin{bmatrix} 0 \\ 0 \\ \frac{1}{4l_c B^2}(F_T^{-1})'(\Psi_0) \end{bmatrix}. \tag{8}$$

Since the operator $A$ is sectorial it generates an analytic semigroup in $X$. We denote by $T(t)$ the semigroup operator on $X$ and the norm and inner product will be denoted by $\|\cdot\|$ and $\langle\cdot,\cdot\rangle$ respectively. In this form the equations can be tackled using the standard optimal control theory in Hilbert spaces.

Observe first that this system is block diagonal. It can be split into two parts: a two dimensional part that describes the evolution of the average flow and the pressue rise, and a part of codimension 2 which describes the evolution of stall. This second part does not depend on the control parameter $\gamma$ and can therefore be integrated separately. In other words, stall does not depend on the control parameter and is therefore uncontrollable. The problem is now reduced to a two dimensional problem.

We seek a feedback control that will minimize the cost functional

$$J(y) = \frac{1}{2}(y_1^2(t_f) + y_2^2(t_f)) + \frac{1}{2}\int_0^{t_f} S(y_1^2(t) + y_2^2(t)) + Ru^2 dt.$$

It is a well known result that the optimal feedback control is given by

$$u = -\frac{1}{R}B^T Q(t) y(t)$$

where the symmetric matrix $Q(t)$ satisfies a matrix Riccati equation, see [7] for more details.

## 5 Basic Control of Design Flow for the Nonlinear Equation

The basic attractor for constant throttle functions has been analyzed completely [4],[6] and one would like to use this knowledge of the asymptotic dynamics when constructiong a control law. However, when $\gamma$ is no longer constant, but a function of the state variables, the components of the basic attractor may change, altering the asymptotic dynamics.

Let us assume for now that we only consider control strategies that move the throttle in an adiabatic fashion. Restricting the control to this class guarantees that the basic attractor is unchanged. The best one can hope to do here is to slide the solution along the basic attractor until it reaches the desired operation point. If stall occurs, then one slides the system along the branch of stable stall cells, by increasing γ, until the saddle-node bifurcation point is reached and stall ceases to exist and the flow converges to design flow. This design flow is achived at a very low pressure rise. In order to increase the pressure rise we decrease γ again until the desired operation point is reached.

Let $\xi(t)$ be a parameterization of a family of stationary solutions in the basic union attractor $\cup_\gamma \mathcal{B}$. Note that since we are working in a rotating frame of reference, $(\eta,t)$, stall solutions will be stationary solutions. Then corresponding to $\xi(t)$ we can find $\bar{\gamma}(t)$ such that for $\gamma = \bar{\gamma}(t)$, $\xi(t)$ is a stationary solution of the equation (1-3). We denote the solution of the system (1-3) as $x(t) = (\Phi, \Psi, \varphi)(t) = \xi(t) + y(t)$ and the control parameter $\gamma = \bar{\gamma}(t) + u(t)$. Let us now linearize this system about the trajectory $\xi(t)$ and write it as

$$\dot{y} = A(t)y + B(t)u - \dot{\xi}(t). \qquad (9)$$

Here $A(t)$ and $B(t)$ depend on time through the trajectory $\xi(t)$. Our goal is to make $y$ as small as possible with very little control effort. In other words, we want to find a regulator for equation (9) which is optimal with respect to the cost function

$$J(u) = \frac{1}{2}\|S_f y(t_f)\|^2 + \frac{1}{2}\int_0^{t_f} \|Sy(t)\|^2 + Ru^2(t)\,dt$$

where $R$ is a constant and $S_f$ and $S$ are symmetric positive definite linear operators. This is a well known problem and it can be solved exactly, see [7].

### 5.1 Construction of the Controller

To make the construction of the controller as simple and intuitive as possible we proceed in the following way. When a disturbance occurs in the system that is large enough so that the system cannot recover without intervention, we change the control parameter to a setting where the only component in the basic attractor is the design flow. This consists of increasing γ to a level $\gamma_1$ so that the throttle characteristic no longer intersects the branch of stall cells. We then wait until the flow is in a small enough neighborhood $U$ of the design flow. This design flow setting is however at a low pressure rise level so to increase the pressure rise we now track a trajectory $\xi_1$ to the desired design flow setting, see Figure 2. As we will prove later, if the state is close enough to the starting point of $\xi_1$ and the cost on the control small enough, this strategy will work for all initial disturbances.

This control construction will still be very close to the original one as the system will settle into stall or surge very fast and then traverse near the basic attractor towards the design flow corresponding to the throttle setting $\gamma_1$.

The linearization of the system about the trajectory $\xi_1$ is exactly that given by equations (7) and (8), except for that now these operators are time-dependent, i.e. instead of $(\Phi_0, \Psi_0)$ we have $(\Phi, \Psi)_{\bar{\gamma}(t)}$. Just like before, there is an uncontrollable subspace of codimension 2, but since we are on the right side of the peak of the characteristic, this space is stable and small disturbances will decay in time. We therefore only need to consider the first two modes which describe the flow in the $(\Phi, \Psi)$-plane and it suffices to know $Q(t)$, $\zeta_1$ and $\zeta_2$. The following theorem is proven in [7].

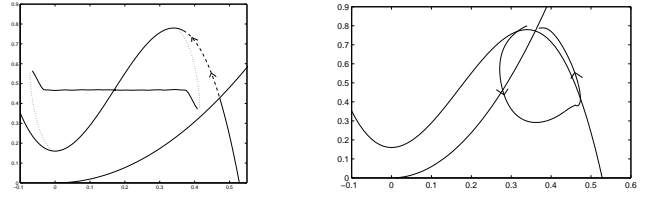

**Figure 2:** Left: The throttle setting $\gamma = \gamma_1$ that defines the start of the trajectory $\xi_1$ which is shown as a dashed line. Right: The $(\Phi, \Psi)$-phase plane for the basic control.

**Theorem 5.1** *There exists an open set U around* $\xi_1(0) = (\varphi_1, \Phi_1, \Psi_1)$, *a constant R and a prevalent set* $Y \subset X$ *such that the control strategy given by the above construction, make all solutions of the viscous Moore-Greitzer equation converge to the desired stable design flow.*

## 6 Numerical Simulations

Here we present some numerical simultions that display how our control performs and its performance is compared with that of the backstepping control given by Banaszuk et al. [1]. For all of the simulations the initial condition is a small disturbance in the average flow and pressure rise, but a large disturbance in the stall direction. The two ordinary differential equations (2) and (3) are solved by a Runge-Kutta routine which is coupled together with a Lax-Wendroff scheme which solves equation (1).

The backstepping control is a much more forceful control that uses more control effort and the pressure rise drops completely when it is used to control stall, see Figure 3. The state has a much smaller excursion in the $(\Phi, \Psi)$-plane with the basic control, see Figures 2, and in particular the pressure rise never drops completely.

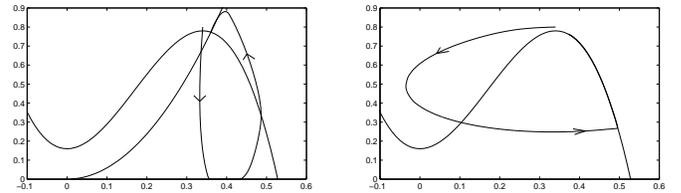

**Figure 3:** Left: The $(\Phi, \Psi)$-phase plane for the backstepping control of stall. Right: The $(\Phi, \Psi)$-plane when controlling a surging compressor with the basic control.

Surge is in general harder to control than stall. It requires more control effort and is a more violent instability. We present here some simulations which show how the two controls handle a surging compressor. Both basic and backstepping controls are saturated to control surge. The control strategy with the least effort that could recover design flow from surge would probably just involve increasing $\gamma$ slightly and then waiting for the system to complete a single surge oscillation. Figure 3 shows the $(\Phi,\Psi)$ phase plane during the transient.

## 7 Conclusion

We defined b-controllability and a-controllability and presented arguments why these would be meaningful definitions of controllability for infinite dimesional nonlinear dynamical systems.

The backstepping control presented by Banaszuk *et al.* was the first attempt at constructing a control strategy for the Moore-Greitzer partial differential equation. The vMG equation, which is a better physical model for the airflow through the compression system, has different asymptotic dynamics than the Moore-Greitzer equation and these asymptotic dynamics have been analyzed by the authors in [4], [5] and [6]. Here and in [7] we go one step further and use the knowledge of the asymptotic dynamics to construct a control strategy that utilizes the dynamics and hence needs considerably less control effort. We believe that this approach of analyzing the basic attractor and using the knowledge of the asymptotic dynamics to construct basic control strategies for b-controllable systems can offer a viable alternative to linearizing high dimensional systems and applying linear optimal control theory to the linearized system.

## 8 Acknowledgements

The authors would like to thank Professor Petar Kokotović for helpful conversations and support. This work was funded by NSF under grant number DMS-9704874 and in part by AFOSR under grant number F49620-95-1-0409. Equipment were funded by NSF under grants number DMS-9628606 and PHY-9601954.